\documentclass[a4paper,10pt]{amsart}

\usepackage{hyperref}

\usepackage{amsmath,amsthm} 
\usepackage{amssymb,mathrsfs} 
\usepackage{bbm}
\usepackage{a4wide} 
\usepackage{physics}
\usepackage{color,subfigure} 
\usepackage{enumerate}
\usepackage[normalem]{ulem}
\usepackage{cancel}
\usepackage[utf8]{inputenc}
\usepackage{comment}

\DeclareMathAlphabet{\mathpzc}{OT1}{pzc}{m}{it}

\newcommand{\cL}{\mathcal{L}}

\newcommand{\cD}{\mathcal{D}}

\renewcommand{\leq}{\leqslant}
\renewcommand{\geq}{\geqslant}

\def\cF {\mathcal{F}}
\def\bR {\mathbbm{R}}
\def\bN {\mathbbm{N}}

\newcommand{\vertiii}[1]{{\left\vert\kern-0.25ex\left\vert\kern-0.25ex\left\vert #1 
    \right\vert\kern-0.25ex\right\vert\kern-0.25ex\right\vert}}

 \newtheorem{Thm}{Theorem}[section]

\newtheorem{Prop}[Thm]{Proposition}
\newtheorem{Cor}[Thm]{Corollary}

\begin{document}

  
\title[]{Loskot-Rudnicki's inequality and General Relative Entropy inequality for Cauchy problems preserving positivity}
\maketitle

\author{\'Etienne Bernard}
\address{CERMICS, \'{E}cole des Ponts, Marne-la-Vallée, France} 

\email{etienne.bernard@enpc.fr}
  

\begin{abstract}
The Generalized Relative Entropy inequality is a ubiquitous property in mathematical models applied in physics or biology. In spite of its importance, it is currently proved on a case-by-case basis in the literature. Here, we show that GRE is actually a generic consequence of Loskot-Rudnicki's inequality that is reminiscent of Jensen's inequality.
\end{abstract}


\maketitle

\section{Introduction}
Generalized Relative Entropy inequality (GRE in short) is a property that often occurs in linear Partial Differential Equations (PDE) preserving the positivity of the initial condition over time. It plays a fundamental role in mathematical models in biology or physics of polymerization (see \cite{perthameTransportEquationsBiology2007,perthameExponentialDecayFragmentation2005,calvoLongTimeAsymptoticsPolymerization2018,jourdainLongTimeAsymptoticsMultiscale2006} and references therein) where it is used to get {\it a priori} estimates, establish a contraction principle or studying long time convergence to a steady state or periodic solutions (e.g. \cite{michelGeneralEntropyEquations2004, michelGeneralRelativeEntropy2005,perthameTransportEquationsBiology2007, bernardCyclicAsymptoticBehaviour2019a}). It has been established for various linear PDE but only on a case-by-case basis. Moreover, in a few papers such as in \cite{bernardCyclicAsymptoticBehaviour2019a}, it remains at a formal level without always avoiding the pitfall of circularity. In \cite{loskotRelativeEntropyStability1991},  Krzysztof Loskot and Ryszard Rudnicki established a Jensen-type inequality for $L^{1}$ spaces and showed that it implies a Csiszár inequality for conservative positive operators. This inequality is in fact a particular case of GRE. The aim of the present paper is to recall in Section \ref{reformulation} Krzysztof Loskot and Ryszard Rudnicki's main results in a slightly more general framework. Then, we show in Section \ref{examples} with some examples how to derive GRE from it in both conservative and non conservative cases. In other words, Loskot-Rudnicki's inequality implies that GRE is a generic propriety for Cauchy problems preserving positivity in a large class of spaces.\\

\section{A Jensen-type inequality in KB spaces}
\label{reformulation}
The results established by Krzysztof Loskot and Ryszard Rudnicki are stated in $L^{1}$ but it is in fact valid for a larger class of spaces that are important in applications such as weighted Lebesgue spaces. So that we will both recall their results and we restate them in a more general framework, the Riesz spaces. For the sake of completeness, we first recall the relevant elements of the theory of Riesz spaces. We refer the reader \cite{schaeferBanachLatticesPositives1974, zaanenIntroductionOperatorTheory1997, p.meyer-niebergBanachLattices1991} for a more comprehensive introduction.\\

 A vector space $E$ over $\bR$ endowed with an order relation $\leq$ is a {\it ordered vector space} if these axioms are satisfied:
 \begin{enumerate}
 \item $x\leq y \Rightarrow x+z\leq y+z\ \mbox{for all }x,y,z\in E,$
 \item $x\leq y \Rightarrow \lambda x\leq \lambda y \ \mbox{for all }x,y\in E\ \mbox{and }\lambda \in\bR_+.$
 \end{enumerate}
 An element $x$ of $E$ is said to be {\it positive} if $0\leq x$ and the subset of positive elements, denoted $E_{+}$, is called the {\it positive cone}.  It is moreover a {\it vector lattice} (also called {\it Riesz space}) if $x \vee y :=\sup\{x,y \}$ and $x \wedge y :=\inf\{x,y \}$ are well-defined in $E$. The absolute value is then defined as $\qty|x| = x^{+} + x^{-}$ with $x^{\pm}:=\qty(\pm x) \vee 0.$ A subset $U$ of $E$ is called {\it solid} if $0\leq x \leq y$ for any $(x,y)\in E\times U$ implies that $x\in U$. An (order) {\it ideal} in $E$ is a subspace that is solid. Let $E,F$ be two ordered space, a linear $U\in \cL\qty(E,F)$ is said to be {\it positive} if it send $U\qty(E_{+})\subseteq F_+.$ A {\it Banach lattice} is a vector lattice with a Banach norm such that $0 \leq x \leq y$ implies $\norm{x} \leq \norm{y}.$ A stochastic operator $U\in \cL\qty(E,F)$ is a positive operator such that $\norm{Ux}_{F}=\norm{x}_{E}$ for any $x\in E_{+}.$\\
 
 A sequence $\qty(x_i)_{i\in \bN}$ of a vector lattice is increasing (respectively decreasing) if $ i\leq j$ implies $x_{i} \leq x_{j}$ (respectively  $x_{j} \leq x_{i}$.)  A sequence is monotone if it is either increasing either decreasing.  A Kantorovich-Banach space (KB-space) is a Banach lattice such that any monotone norm bounded sequence is convergent.  It is known that $E$ is a KB-space if and only if it is weakly sequentially complete (Theorem 2.5.6 p. 104 in  \cite{p.meyer-niebergBanachLattices1991}). In particular, reflexive Banach lattices and Banach lattice with a p-additive norm is a KB-space (Corollary 2.413 p. 93 in  \cite{p.meyer-niebergBanachLattices1991}), including therefore every $L^{p}$ spaces for $1\leq p < +\infty$.\\

 Let $\qty(X,\Sigma,\mu)$ be a complete $\sigma-$finite measure space, then the set $L^{0}\qty(X,\Sigma,\mu)$ of all $\Sigma-$measurable $\mu-$almost everywhere finite real valued functions modulo $\mu$-null functions endowed with the pointwise order ($f \leq g$ if and only if $f(x)\leq g(x)$ $\mu-$a.e.) is a Riesz space (see \cite{zaanenIntroductionOperatorTheory1997} p. 12). To reduce the amount of notation, we drop henceforth $\qty(X,\Sigma,\mu)$ and we write only $L^{0}$ for $L^{0}\qty(X,\Sigma,\mu).$  The ideals of $L^{0}$ are called in literature {\it function spaces}, see e.g. \cite{abramovichInvitationOperatorTheory2002} p.~194. The spaces of (bounded, vanishing at infinity, and so on) continuous functions are not ideals of $L^{0}$ and therefore are {\it not} function spaces. We call {\it function} KB-space a function space that is also a KB space. Indeed, the Lebesgue space $L^p\equiv L^{p}\qty(X,\Sigma,\mu)$ ($1\leq p < \infty$) are function KB-spaces.\\

 Henceforth, let $E$ be a function KB-space. For any positive linear operator $U\in \cL\left(E\right)$, we denote $\tilde{U}$ its extension to $L^{0}$ by
$$
\tilde{U}f := \sup \qty{Ug: g \in E, g \leq f },\ f\in L^0.
$$ 
Let $\eta$ be a continuous convex function defined on $[0,\infty)$,  the authors of \cite{loskotRelativeEntropyStability1991} introduce the auxiliary function $\phi_{\eta}: [0,\infty)\times [0\infty) \rightarrow \overline{\bR}$
$$
\phi_{\eta}\qty(u,v):=
\begin{cases}
v \eta\qty(\frac{u}{v}), &v>0, u\geq 0,\\
0, &v=0, u=0,\\
u\eta^{'}(\infty), &v=0, u>0,
\end{cases}
$$
with $\eta^{'}(\infty) = \lim_{x\rightarrow +\infty}\frac{\eta(x)}{x}$  and proved a very interesting inequality (Proposition 2.2 of \cite{loskotRelativeEntropyStability1991}): 
\begin{Prop}[Loskot-Rudnicki's inequality]
\label{Jensen}
Let $E$ be a KB function space and let $U\in \cL_{+}\qty(E)$. Then 
$$
\phi_{\eta}\qty(U f , U g) \leq \tilde{U} \phi_{\eta} \qty(f , g). 
$$
\end{Prop}
The authors have established it for $L^{1}$ but a close inspection of their proof and the remark that they have used the convergence monotone theorem and the property of convex as supremum of affine functions show that it can be easily extended to the case of function KB-spaces. With this inequality, they establish Csiszár's inequality for every stochastic operator U from $L^{1}\qty(\mu_1)$ to $L^1\qty(\mu_2)$ (Theorem 2.1 in \cite{loskotRelativeEntropyStability1991}): 
\begin{Cor}
\label{LoskotGRE}
Let $U: L^{1}\qty(\mu_1) \rightarrow L^{1}\qty(\mu_1)$ being a stochastic operator and let $\eta:\bR_{+} \rightarrow \overline{\bR}$ being a convex function. Let denote
$$
H_{\eta,i}\qty(f | g):= \int \phi_{\eta}\qty(f,g) d\mu_i,\ \forall f,g \in L^{1}\qty(\mu_i).
$$
Then we have
$$
H_{\eta,2}\qty(Uf | U g) \leq H_{\eta,1}\qty(f|g)\ \forall f,g \in L^{1}\qty(\mu_1).
$$
\end{Cor}
Corollary \ref{LoskotGRE} implies GRE in the conservative case as we will see in the first example of Section \ref{examples}. But Loskot-Rudnicki's inequality is more important as it implies some version of GRE in non-conservative cases, see the second example in Section \ref{examples}.
 
 \section{Examples}
 \label{examples}
  \subsection{The growth models}
  \label{growthmodel}
As first example, we consider one of the first problems that motivated the introduction of GRE, i.e. from \cite{michelGeneralRelativeEntropy2005} a growth model which "take the form of a mass preserving fragmentation equation with a drift term" applied to the time dynamic of the population of particles/cells/individuals:
 \begin{equation}
 \label{growtheq}
 \begin{cases}
\pdv{t}n + \cD_0 n = \cF\ \mbox{on }(0, \infty)\times(0,\infty)\\
\mbox{boundary condition in }x=0\\
n(t=0,x) = n_{0}(x)
 \end{cases}
 \end{equation}
 where $\cF$ is a mass conservative fragmentation operator 
 $$
 \qty(\cF n)(t,x) = \int_{0}^{\infty} b(t,y,x) n(t,y)dy - n(t,x) B(t,x)
 $$
 and $\cD_0$ the drift term with velocity $v(x)\geq 0:$
 $$
 \qty(\cD_{0} n)(t,x) = \pdv{x}\qty(v(x)n(t,x)) + w(t,x)n(t,x).
 $$
 We refer the reader to \cite{michelGeneralRelativeEntropy2005} for an explanation of the biological motivations behind the operators introduced above. Let introduce the associated dual equation:
 $$
 -\pdv{t}\psi(t,x) + \cD_{0}\psi = \cF^{*}\psi
 $$
We assume that the coefficients are such that (Theorem 2.2 in \cite{michelGeneralRelativeEntropy2005}) there exists $\psi >0$ being a solution to the dual equation with initial condition $\psi(0,) = \psi_{0}$, then for any initial datum $n_{0}$ such that $n_{0}\psi_{0} \in L^{1}\qty(0,\infty)$, there exists a (unique) solution to equation \ref{growtheq} such that 
\begin{equation}
\label{conservative}
\int_{0}^{\infty} n(t,x) \psi(t,x) dx = \int_{0}^{\infty} n_{0}\psi_{0}dx.
\end{equation}
Let define the operator $T_{t}: n_{0} \in L^{1}\qty(\psi_{0}(x)dx) \mapsto n(t, \cdot) \in L^{1}\qty(\psi(t, x)dx).$ It is a positive operator and moreover stochastic by $\qty(\ref{conservative})$. Let be $H$ a convex function defined on $[0,+\infty).$ Then by Corollary \ref{LoskotGRE}, we recover the decreasing of GRE ((1.3) in \cite{michelGeneralRelativeEntropy2005}):
$$
\int_{0}^{\infty}\psi(t,x) \phi\qty(T_t n_{0}, T_{t} m_{0})dx \leq \int_{0}^{\infty} \psi_{0} \phi\qty(n_0, m_0)(x)dx,\ \forall t\geq0.
$$
\subsection{An example of non-conservative case: the transport equation in a domain with absorbing boundary}
 We show here with an example how we can derive GRE from Proposition \ref{Jensen} in non-conservative cases. We consider the transport problem in a bounded domain with absorbing boundary. Let $\Omega \subset \bR^{n}$ with smooth boundary $\Gamma$ be the spatial domain, and let $V\subset \bR^{n}$ the velocity domain. We denote for any $x\in \Gamma,$ $n_x$ the outer normal vector. Consider $f\equiv f(t,x,v)$ being solution of the following transport problem:
 \begin{equation}
 \label{transport}
 \begin{cases}
 \pdv{t} f + v\cdot\nabla_{x} f + \sigma\qty(f - Kf)=0, \\
 f_{\left|n_x\cdot v < 0\right.} =0, \\
 f(0,x,v) = f_{0}(x,v),
 \end{cases}
 \end{equation}
 with $K\in \cL\qty(L^{1}\qty(V)).$ It can model for instance the time evolution of the motion of neutrons in an absorbing and scattering homogeneous medium, with  $f$ being the density of population of particles. The boundary condition means that any particle hitting $\Gamma$ exits $\Omega$ forever (see chapter 17 in \cite{a.batkaiPositiveOperatorSemigroups2017}).  It is known that it generates a positive $C_{0}$-semigroup $\qty(T_t)_{t\geq0}$ in $L^{1}\qty(\Omega \times V)$ (see \cite{arlottiNewApproachTransport2009}) with a loss of mass:
 \begin{equation}
 \label{conv}
 \iint_{\Omega \times V}\qty(T_{t} f_{0})(x,v) dxdv \leq  \iint_{\Omega \times V} f_{0}(x,v) dxdv.
 \end{equation}
Let consider $H: \bR_{+} \rightarrow \overline{\bR}$ a convex function bounded from below by a constant, and let $f_{0}, g_{0} \in L^{p}\qty(\Omega \times V)$ with $g_{0} >0$ a.e. and such that $g_0 H\qty(\frac{f_0}{g_0})\in L^{1}\qty(\Omega \times V).$ Let $f = T_t f_0$ and $g=T_t g_0.$ To avoid technicalities, we assume that $g>0$ a.e.. By Proposition \ref{Jensen}, we have 
 $$
g H\qty(\frac{f}{g})\leq T_t \qty(g_0 H\qty(\frac{f_0}{g_0})).
 $$
 Let $C$ a constant such that $H+C \geq 0$ then the inequality implies that
 $$
 0\leq Cg + g H\qty(\frac{f}{g}) \leq Cg + T_t\qty(g_0 H\qty(\frac{f_0}{g_0}))
 $$
 As $Cg + T_t \qty(g_0 H\qty(\frac{f_0}{g_0})) \in L^{1}\qty(\Omega \times V)$ that is an ideal of $L^0$, that implies that $Cg + g H\qty(\frac{f}{g}) \in L^{1}\qty(\Omega \times V)$ and thus $g H\qty(\frac{f}{g}) \in L^{1}.$ Finally we get by $(\ref{conv})$
 $$
 \iint_{\Omega \times V}g H\qty(\frac{f}{g}) \leq \iint_{\Omega \times V} g_0 H\qty(\frac{f_0}{g_0})
 $$
 Therefore, we have the following proposition:
 \begin{Prop}
 Let $f$ and $g$ being mild solutions of $\ref{transport}$ with respective initial conditions $f_0, g_0$. Assume that $g >0$ a.e. for all $t\geq0$, then if $g_0 H\qty(\frac{f_0}{g_0}) \in L^{1}\qty(\Omega \times V)$, then $\forall t>0$, $g H\qty(\frac{f}{g}) \in L^{1}\qty(\Omega \times V)$ and the relative entropy $t\mapsto \iint_{\Omega \times V} T_t g_0 H\qty(\frac{T_t f_0}{T_t g_0})$ is nonincreasing over time.
  \end{Prop}
Notice that in both case, the semigroup property and Proposition \ref{Jensen} imply GRE. However, the regularity of the semigroup do not intervene, therefore we can derive GRE for discrete semigroups in the same way than above. That is an advantage in the study of discrete-time dynamical systems coming from biological or physical models, or more widely dynamical systems with low regularity. 
%

\section*{Acknowledgements}

 The author is grateful to Benoît Perthame for bringing to his attention the subject of Generalized Relative Entropy inequality and for useful discussions. The author also thank Pierre Gabriel for useful discussions and Tony Lelièvre for also useful discussions and for having pinpointed \cite{jourdainLongTimeAsymptoticsMultiscale2006}. Eventually, he also thank the anonymous reviewers for their careful reading and their many comments and suggestions, especially pinpointing \cite{loskotRelativeEntropyStability1991} while recalling the concept of KB-spaces. I also thank them for suggesting some very useful references.  All their advice has greatly shortened and improved the manuscript.

\bibliographystyle{plain}
\bibliography{Maths.bib}

\end{document}